\documentclass[a4paper,leqno,11pt,final]{article}
\usepackage{latexsym}
\usepackage{amsthm,amsmath,amssymb}

\usepackage[dvips]{graphicx, color}

\newtheorem{thm}{Theorem}

\newtheorem{fact}[]{Fact}

\newtheorem{prop}[thm]{Proposition}

\newtheorem{conj}[thm]{Conjecture}

\newcommand{\Aut}{\mathrm{Aut}}

\begin{document}
\title{\bf {\Large{General Uniformity of Zeta Functions}}}
\author{\bf Lin WENG}
\date{\bf }
\maketitle
\hskip 5.0cm --- In memory of Prof. Dr. F. Hirzebruch
\vskip 0.30cm
\noindent
{\bf Abstract.} {\footnotesize Using analytic torsion associated to stable bundles, we introduce zeta functions for compact Riemann surfaces.
To justify the well-definedness, we analyze the degenerations of analytic torsions at the boundaries of the moduli spaces, the singularities of analytic torsions  at Brill-Noether loci, and  the asymptotic behaviors of analytic torsions with respect to the degree.
These new yet intrinsic zetas, both abelian and non-abelian, are expected to play key roles to understand global analysis and geometry of Riemann surfaces, such as the Tamagawa number conjecture for Riemann surfaces, searched by Atiyah-Bott, and the volumes formula of moduli spaces of Witten. Relating to this, in our theory on  special uniformity of zetas, we will first construct a symmetric zetas based on abelian zetas and group symmetries, then  conjecture that
our non-abelian zetas coincide with these later zetas with symmetries. 
All this, together with that for zetas of number fields and function fields, then consists of our theory of general uniformity of zetas.}
\tableofcontents

\section{Zeta Functions for Riemann surfaces}
\subsection{Regularized Integrations}
Let $M$ be a compact Riemann surface of genus $g>1$ with fundamental group $\pi_1(M)$. Fixed a normalized smooth volume form $\omega$ on $M$ so that $\int_M\omega=1$. For a line bundle $L$ of degree $d$, by definition, an $\omega$-admissible metric $h$ on $L$ is a hermitian metric  on $L$ such that $c_1(L,h)=d\cdot\omega$. One checks that $\omega$-admissible metrics always exist and for a fixed line bundle they are parametrized by positive reals. For our purpose, fix a point $P_0\in M$ and use the normalized Green function ([L]) to define the $\omega$-admissible metric $h_A$ on the line bundle $A=A_M$ corresponding to the invertible sheaf $\mathcal O_M(P_0)$. (From now on, we will not make distinctions between bundles and locally free sheaves.) Moreover, fix a conformal metric $\tau=\tau_M$ on $M$ such that its induced metric on the canonical line bundle $K_M$ is $\omega$-admissible and whole volume is given by $2\pi(2g-2)$.

Let $\mathcal M_{M,r}(d)$ be the moduli space of stable bundles of rank $r$ and degree $d$ on $M$. Then we have natural isomorphisms, for all $m\in\mathbb Z$,
$$\mathcal M_{M,r}(0)\simeq \mathcal M_{M,r}(mr),\qquad V\mapsto V\otimes A^{\otimes m}.$$
By a result of Narasimhan-Seshadri [NS], points $V$ of $\mathcal M_{M,r}(0)$ are in one-to-one correspondence to irreducible unitary representations $\rho$ of $\pi_1(M)$. (Here for simplicity, we assume that $g\geq 2$.) This will be indicated as $V=V_\rho$. Use the uniformization, denote the induced hermitian metric by $h_\rho$ on $V_\rho$ from the standard one on $\mathbb C^r$. As such, then for any $V\in \mathcal M_{M,r}(rm)$, we have the associated canonical metric $h_V=h_\rho\otimes h_A^{\otimes m}$ on $V=V_\rho\otimes A^{\otimes m}$ induced from the canonical metric $h_\rho$ on $V_\rho$ and the canonical metric $h_A$ on $A$ constructed above. For simplicity, denote such a metrized bundle (resp. Riemann surface $M$) by $\overline V$ (resp. $\overline M$).

Denote by $\tau(M,V)$ the Ray-Singer analytic torsion associated to the hermitian vector bundle $\overline V$ over $\overline M$. Then the determinant line bundle $\lambda$ on $\mathcal M_{M,r}(rm)$, 
whose fiber at $V\in \mathcal M_{M,r}(rm)$
is given by $\det H^0(M,V)\otimes\det H^1(M,V)^{\otimes -1}$, together with the associated Quillen metric $h_Q(V,M)$, induced from that of $\overline V/\overline M$, is a definite metrized polarization on $\mathcal M_{M,r}(rm)$. As above, denote the resulting metrized line bundle by $\overline \lambda$.

Let $\mathcal V_m$ be the universal Poincare vector bundle on $\mathcal M_{M,r}(rm)$. Then with respect to the projection $q:M\times \mathcal M_{M,r}(rm)\to \mathcal M_{M,r}(rm)$, {\it assume} that,
from the short exact sequence of
coherent sheaves $$0\to \mathcal V_m\to \mathcal V_m\otimes A^{\otimes d}\to \mathcal V_m\otimes A^{\otimes d}/\mathcal V_m\to 0,$$ we obtain the exact sequence  of {\bf vector bundles} for 
higher direct images
$$0\to q_*\mathcal V_m\to q_*\Big(\mathcal V_m\otimes 
A^{\otimes d}\Big)\buildrel \gamma\over\to q_*\Big(\mathcal V_m\otimes A^{\otimes d}/\mathcal V_m\Big)\to R^1q_*\mathcal V_m\to 0.$$ Denote by $W_{M,r}^{\geq i}(d)$ the determinantal variety associated to $\gamma$. Then one checks that
the support of $W_{M,r}^{\geq i}(d)$ coincides with the so-called Brill-Noether locus consisting of
these whose $h^0$ are at least $i$. (See e.g., [ACGH, p.176].)
Indeed, $W_{M,r}^{\geq i}(d)$
are normal subvarieties of $\mathcal M_{M,r}(d)$ and $T(V,M)$ is a smooth function on $$W_{M,r}^{i}(d)=W_{M,r}^{\geq i}(d)\,\backslash\,W_{M,r}^{\geq (i+1)}(d)$$ when $d\in r\mathbb Z$. 
By an abuse of notation, denote by $d\mu$ the volume forms on $W_{M,r}^{\geq i}(d)$ induced from the polarization $\lambda$. 

With all this, we are ready to introduce the regularized integration by
$$\begin{aligned}\int_{\mathcal M_{M,r}(rm)}^{\#}&\Big(e^{\tau(M,V)}-1\Big)\big(e^{-s}\big)^{\chi(M,V)}d\mu\\
:=\sum_{i=0}^\infty &\int_{W_{M,r}^{i}(rm)}\Big(e^{\tau(M,V)}\Big)\big(e^{-s}\big)^{\chi(M,V)}d\mu-\int_{\mathcal M_{M,r}(rm)}\big(e^{-s}\big)^{\chi(M,V)}d\mu.\end{aligned}$$

\subsection{Zetas for Riemann Surfaces}
Let $M$ be a compact Riemann surface and fix $r\in\mathbb N$. 
We define {\it the rank $r$ zeta function for $M$} by
$$\widehat\zeta_{M,r}(s):=\sum_{m=-\infty}^\infty \int_{\mathcal M_{M,r}(rm)}^{\#}\Big(e^{\tau(M,V)}-1\Big)\big(e^{-s}\big)^{\chi(M,V)}d\mu,\qquad\mathrm{Re}(s)>0.$$
The main result of this paper is the following
\vskip 0.30cm
\begin{thm} ({\bf Zeta Facts}) (i) $\widehat\zeta_{M,r}(s)$ is well defined and admits a (unique) meromorphic continuation to the whole complex $s$-plane.

\noindent
(ii) ({\bf Functional Equation})
$$\widehat\zeta_{M,r}(-s)=\widehat\zeta_{M,r}(s).$$

\noindent
(iii) ({\bf Singularities)} $\widehat\zeta_{M,r}(s)$ admits only one singularity, namely, a simple pole at $s=0$ with the
residue $$\mathrm{Res}_{s=0}\widehat\zeta_{M,r}(s)=\mathrm{Vol}\Big(\mathcal M_{M,r}(0)\Big)$$ where 
$\mathrm{Vol}\Big(\mathcal M_{M,r}(0)\Big)$ denotes the volume of the moduli space
$\mathcal M_{M,r}(0)$ with respect to the volume $d\mu$.
\end{thm}

\subsection{Zeta Facts I: Formal Aspect}
Before justifying the convergence of our regularized integrations appeared in the definition of new zeta functions, let us formally establish the functional equation and find out the singularities of these zetas and hence calculate the associated residues when applicable.

For this purpose, we first recall some basic facts about analytic torsions. Let $\overline V/\overline M$ be a metrized vector bundle over a metrized Riemann surface. Denote the associated Laplacian by $D_V$ on the associated space of $L^2$  sections $L^2(\overline M,\overline V)$ of $V$ on $M$. From the Fredholm theory, the spectrum of $D_V$ is a purely discrete sequence
$$0\leq \lambda_1\leq\lambda_2\leq\lambda_3\leq\dots,\qquad\lambda_n\sim\frac{1}{r(g-1)}n$$ with corresponding eigenfunctions $\{e_n(z,V)\}$ forming a complete orthonormal basis.
Accordingly, we can define for $\mathrm{Re}(\lambda)>0$ and $\mathrm{Re}(s)>1$ the spectrum zeta function $$\zeta_\lambda(s,V):=\mathrm{Tr}(D_L+\lambda)^{-s}=\sum_{n=1}^\infty\frac{1}{(\lambda+\lambda_n)^s}
$$ and more generally for any $c\geq 0$,
$$\zeta_\lambda^c(s,V):=\sum_{\lambda>c}\frac{1}{(\lambda+\lambda_n)^s}.$$

We have the following

\begin{thm} ([RS]) (i) For fixed $c\geq 0$, $\zeta_\lambda^c(s,V)$ has an analytic continuation to the half plane $\mathrm{Re}(\lambda)\geq c$ and a meromorphic continuation to the whole $s$-plane with only a simple pole at $s=1$ with residue $r(g-1)$.

\noindent
(ii) For $\mathrm{Re}(\lambda)>\lambda_*$ the smallest non-zero eigenvalue of $D_V$,
$$\zeta_\lambda^0(s,V)=\zeta_\lambda(s,V)-h^0(M,V)\lambda^{-s}$$ has an analytic continuation through the $s$-plane with
$$\zeta_\lambda^0(0,V)+h^0(M,V)=-\Big(\lambda+\frac{1}{3}\Big)n(g-1)+\frac{1}{2}d(V).$$

\noindent
(iii) ({\bf Duality}) Equipped the dual bundle $V^\vee$ with the dual metric, then
$$\zeta_\lambda^c(s,V)=\zeta_\lambda^c(s,K_M\otimes V^\vee)$$
\end{thm}

Based on this, following  Ray-Singer ([RS]), define the {\it analytic torsion} for $\overline V(/\overline M)$ by $$T(V):=T(\overline M;\overline V):=e^{\tau(V)}$$ with
$$\tau(V):=\tau(M,V):=\tau(\overline M;\overline V):=\frac{d}{ds}\zeta_\lambda^0(s,V)|_{s=\lambda=0}.$$
It is well known then that $T(V)$ may be viewed as a regularized determinant of the Laplacian $D_V$. That is to say, formally, we have
$$T(V)=\mathrm{det}'(D_V)={\prod}'_{n:\,\lambda_n>0}\lambda_n.$$

\begin{prop} ({\bf Functional Equation}) Assuming the convergence,
$$\widehat\zeta_{M,r}(-s)=\widehat\zeta_{M,r}(s).$$
\end{prop}

To establish it, motivated by our works for number fields ([W]), let us divide the summation in the definition of our zetas into the following groups:
$$\begin{aligned} \widehat\zeta_{M,r}(s):=&\sum_{m\in\mathbb Z}\int_{\mathcal M_{M,r}(rm)}^{\#}\Big(e^{\tau(M,V)}-1\Big)(e^{-s})^{\chi(M,V)}d\mu\\
=&I(s)+II(s)+III(s)-IV(s)\end{aligned}$$
with $$\begin{aligned}
I(s)=&\sum_{m=0}^{2g-2}\int_{\mathcal M_{M,r}(rm)}^{\#}e^{\tau(M,V)}(e^{-s})^{\chi(M,V)}d\mu,\\
II(s)=&\sum_{m<0}\int_{\mathcal M_{M,r}(rm)}\Big(e^{\tau(M,V)}\Big)(e^{-s})^{\chi(M,V)}d\mu,\\
III(s)=&\sum_{m>2g-2}\int_{\mathcal M_{M,r}(rm)}\Big(e^{\tau(M,V)}-1\Big)(e^{-s})^{\chi(M,V)}d\mu,\\
IV(s)=&\sum_{m\leq 2g-2}\int_{\mathcal M_{M,r}(rm)}(e^{-s})^{\chi(M,V)}d\mu
\end{aligned}
$$
First for $I$, we divide it further as follows:
$$\begin{aligned}
I(s)=&\sum_{m=0}^{2g-2}\int_{\mathcal M_{M,r}(rm)}^{\#}e^{\tau(M,V)}(e^{-s})^{\chi(M,V)}d\mu\\
=&\Big(\sum_{m=0}^{(g-1)-1}+\sum_{m=g-1}
+\sum_{m=g}^{2g-2}\Big)\int_{\mathcal M_{M,r}(rm)}^{\#}e^{\tau(M,V)}(e^{-s})^{\chi(M,V)}d\mu\\
=&\sum_{m=0}^{(g-1)-1}\int_{\mathcal M_{M,r}(rm)}^{\#}e^{\tau(M,V)}(e^{-s})^{\chi(M,V)}d\mu\\
&\qquad+\sum_{m=(g-1)-1}^{0}\int_{V\in\mathcal M_{M,r}(rm)}^{\#}e^{\tau(M,K_M\otimes V^\vee)}(e^{-s})^{\chi(M,K_M\otimes V^\vee)}d\mu\\
&+\sum_{m=g-1}\int_{\mathcal M_{M,r}(rm)}^{\#}e^{\tau(M,V)}(e^{-s})^{0}d\mu\\
=&\sum_{m=0}^{(g-1)-1}\int_{\mathcal M_{M,r}(rm)}^{\#}e^{\tau(M,V)}\Big((e^{-s})^{\chi(M,V)}+(e^{s})^{\chi(M,V)}\Big)d\mu\\
&+\sum_{m=g-1}\int_{\mathcal M_{M,r}(rm)}^{\#}e^{\tau(M,V)}d\mu\\\end{aligned}$$
by duality, i.e., Thm 2.iii) for analytic torsions and the ordinary duality on cohomology groups.
Thus clearly, $I(-s)=I(s).$

Next, for $II(s)$, using the duality, we have
$$\begin{aligned}
II(s)
=&\sum_{m<0}\int_{V\in\mathcal M_{M,r}(rm)}\Big(e^{\tau(M,K_M\otimes V^\vee)}\Big)(e^{-s})^{-\chi(X,K_M\otimes V^\vee)}d\mu\\
=&\sum_{m>2g-2}\int_{\mathcal M_{M,r}(rm)}\Big(e^{\tau(M,V)}\Big)(e^{s})^{\chi(M,V)}d\mu\\
=&III(-s)+V(s)\end{aligned}$$
where
$$V(s)=\sum_{m> 2g-2}\int_{\mathcal M_{M,r}(rm)}(e^{s})^{\chi(M,V)}d\mu.$$
Finally, let us calculate $IV(s)$ and $V(s)$. By definition and the Riemann-Roch, we have
$$\begin{aligned}IV(s)=&\sum_{m\leq 2g-2}\int_{\mathcal M_{M,r}(rm)}(e^{-s})^{\chi(M,V)}d\mu
=\sum_{m\leq 2g-2}\int_{\mathcal M_{M,r}(rm)}(e^{-s})^{r[m-(g-1)]}d\mu\\
=&\sum_{m\leq 2g-2}(e^{-s})^{r[m-(g-1)]}\cdot\mathrm{Vol}\Big({\mathcal M_{M,r}(0)}\Big)\\
=&\frac{(e^{rs})^{(1-g)}}{1-e^{rs}}\cdot\mathrm{Vol}\Big({\mathcal M_{M,r}(0)}\Big)=
\frac{(e^{-rs})^{g}}{e^{-rs}-1}\cdot\mathrm{Vol}\Big({\mathcal M_{M,r}(0)}\Big)
\end{aligned}$$ and similarly
$$\begin{aligned}V(s)=&\sum_{m> 2g-2}\int_{\mathcal M_{M,r}(rm)}(e^{s})^{\chi(M,V)}d\mu\\
=&\sum_{m>2g-2}(e^{s})^{r[m-(g-1)]}\cdot\mathrm{Vol}\Big({\mathcal M_{M,r}(0)}\Big)\\
=&\frac{(e^{rs})^{g}}{1-e^{rs}}\cdot\mathrm{Vol}\Big({\mathcal M_{M,r}(0)}\Big)
\end{aligned}$$
Clearly, $$V(-s)=-VI(s).$$
This then formally completes the proof of the functional equation.
Consequently, if $III(s)$ admits an analytic continuation to the whole $s$-plane, a fact which will be established in the next subsection, we then conclude the following
\begin{prop} ({\bf Singularity and Residue}) Assume the convergence. Then
$\widehat\zeta_{M,r}(s)$ admits only a single singularity, namely, a simple pole at $s=0$
with residue $$\mathrm{Res}_{s=0}\widehat\zeta_{M,r}(s)=\mathrm{Vol}\Big({\mathcal M_{M,r}(0)}\Big)
.$$
\end{prop}

\noindent
{\it Remarks.} This definition of our zetas for Riemann surfaces
are motivated by our works on zetas for numbers and curves over finite fields ([W]).
To understand our current definition, we mention the follows:

\noindent
(i) Instead of working on $\mathcal M_{M,r}(d)$ for all $d$, we take only these which are multiples of the rank $r$, namely, $d\in r\mathbb Z$. This is due to the fact that over finite fields, such and only such a pure selection would guarantee the Riemann Hypothsis.

\noindent
(ii) The summation on all $d\in r\mathbb Z$, instead of  $d\in r\mathbb Z_{\geq 0}$, is motivated by our zetas for number field: after all analytic torsions do not satisfy the vanishing but satisfy the duality.

\subsection{Zeta Facts II: Analytic Aspect}
To establish the zeta facts for our zeta functions, two types of convergences should be justified properly. Namely, the one for  regularized integrations over the moduli spaces $\mathcal M_{M,r}(rm)$ for a fixed $m$, and the other for the infinite sum on $m, m>2g-2$ appeared in $III(s)$. As we will see below,
these two are very different in nature: Technically, for the first type, we need to see how the analytic torsions $T(V)$ degenerate when $h^0$ jump; while for the second, we need to understand how the analytic torsions $T(V\otimes A^{\otimes m})$ behave when 
$m\to\infty$.

\subsubsection{Degenerations of Analytic Torsions}
In this subsection, we will establish the convergence of the regularized integration
$$\int_{\mathcal M_{M,r}(rm)}^{\#}\Big(e^{\tau(M,V)}-1\Big)\big(e^{-s}\big)^{\chi(M,V)}d\mu$$ for each fixed $m$. There are a few issues here:
First, $\mathcal M_{M,r}(rm)$ are not compact; second, $W_{M,r}^{\geq i}(rm)$ are not smooth; and finally, $T(V)$ are not smooth on $W_{M,r}^{\geq i}(rm)$.
\vskip 0.30cm
\noindent
{\bf A. Non-Compactness}
This is not really that serious, thanks to the classical works done.
In fact, using Mumford's GIT, we now have a natural compactification
$\overline {\mathcal M_{M,r}(rm)}$ in terms of Seshadri's equivalences of semi-stable bundles. Denote by $\partial \Big(\mathcal M_{M,r}(rm)\Big):=\overline{\mathcal M_{M,r}(rm)}\backslash \mathcal M_{M,r}(rm)$ the associated boundary. It is well-known that this boundary is of much higher co-dimension, and hence does not cause any serious trouble for the integration.
\vskip 0.30cm
\noindent
{\bf B. Geometric Singularities}

Note that for any open subset $U\subset\mathcal M_{M,r}(rm)$,
if $h^0$ is a constant on $U$, then $T(V)$ is smooth on $U$. Thus we only need to consider the case when $0\leq m\leq g-1$ with the duality 
for analytic torisons in mind.

From the structure of determinantal varieties, it is well-known that the singularities of $W_{M,r}^{\geq i}(rm)$ is contained in $W_{M,r}^{\geq (i+1)}(rm)$, which has codimension at least 2 except in the case when $m=g-1$. Thus, in pour discussion,
we will use the most complicated level, i.e., $m=g-1$, to show how the convergence can be established. For other levels, which are much simpler, a consideration using the following insertion formula for analytic torsions ([AGBMNV]) is sufficient to complete the argument. (For unknown notations, please consult [F].)

\begin{thm} (see e.g. [F, Thm 4.13]) Let $L$ be a line bundle of degree $d\geq g$ with an admissible metric $h$. Then for all stable bundles $V=V_\rho$ with uniformizing metric, $h^1(M,V_\rho\otimes L)=0$ and for any points $x_1,\dots, x_{N:=d+1-g}\in M$:
$$T(V_\rho\otimes L)=\varepsilon_d(M)\frac{T(V_\rho\otimes L(-\sum_{i=1}^NP_i))}{\det(B(x_i,\overline x_j;V_\rho\otimes L))}\cdot\frac{\prod_{i<j}^NP(x_i,x_j)^{2r}}{\prod_{i=1}^mh(x_i)^r}$$ where $\varepsilon_d(\rho)$ is a constant depending only on $d$ and $M$.
\end{thm}
\vskip 0.30cm
\noindent
{\bf C. Analytic Singularities}

So from now on, we concentrate on the level $g-1$. For this, we  recall some facts on both abelian and  non-abelian theta functions. The explosion here follows closely that of Fay [F] (please check the meaning of the unknown notations below in [F] as well).

\begin{thm} ({\bf Theta Functions} [F, Thm 1.6])
Let $L$ be a fixed line bundle with $h^0(\chi(s)\otimes L)$ constant for $s$ in some neighborhood $V$ containing $0\in\mathbb C^N$;  choose 
$\{\omega_i\},\,\{\omega_i^*\}$ bases for $H^0(\chi(s)\otimes L)$, $H^0(\chi^*(s)\otimes KL^{-1})$ with $\{M^{-1}(z,s)\omega_i\},\,\{{^t}M(z,s)\omega_i^*\}$ holomorphic in $s\in V$. Then 

\noindent
(i) for $s\in V$, $$T(\chi\otimes L)=U(\chi(s))|f(s)|^2\det(\langle \omega_i,\omega_j\rangle_{\chi(s)\otimes L})\,\det(\langle \omega_i^*,\omega_j^*\rangle_{\chi^*(s)\otimes KL^{-1}})$$ where $f(s)$ is a holomorphic function on $V$ depending on $L$, the fixed potential $U$, the bases $\{\omega_i\}$, $\{\omega_j^*\}$, and the metrics $h, I\otimes h$ and $\rho$. In particular,

\noindent
(ii) in a neighborhood of any point $\chi(0)$ where $h^0(\chi(s)\otimes \Delta)=0$ with  $\Delta$  a Riemann divisor class satisfying $h^0(M,\Delta)=0$,
$$T(\chi(s)\otimes L)=c_0^l(h,\rho)\,U(\chi(s))\,|\theta(\chi(s))|^2$$ with $\theta(\chi(s))=\theta(s)$ holomorphic in $s$ and independent of the metrics $h,\rho$.
\end{thm}

\begin{thm} ({\bf Vanishing of Non-Abelian Theta},  [F, Prop 4.7, Thm 4.8])
 There exists a holomorphic section $\theta_r$ of the determinant line bundle $\lambda$ on $\mathcal M_{M,r}(0)$ such that 

\noindent
(i) $\theta_r(V_\rho)=0$ if and only if $h^0(M,\mathrm{End}V_\rho\otimes\Delta)>0$;

\noindent
(ii) $\frac{\theta_r(E_\rho)\theta(\det E_\rho)^2}{\theta(\rho)^{2r}}$ is a meromorphic function on $\mathcal M_{M,r}(0)$;

\noindent
(iii) As a bundle on $\mathcal M_{M,r}(0)$, $\lambda\simeq K_{\mathcal M_{M,r}(0)}^{\vee}$, the dual of the canonical line bundle of $\mathcal M_{M,r}(0)$.

\noindent
(iv) The section $\theta_r$ vanishes to order $n$ at any representation $V_\rho\in \mathcal M_{M,r}(0)$
with $h^0(V_\rho\otimes\Delta)=n$. The tangent cone to $(\theta_r)$ at $V_\rho$ is the sub variety of $\overline w=\sum_{i=1}^{\mathrm{dim}\mathcal M_{M,r}(0)}s_i\overline w_i(z,\mathrm{End}V_\rho)\in\overline{H^0(M,K_M\otimes\mathrm{End}V_\rho)}$ given by
$$\det_{1\leq i,j\leq n}\Big(\int_M{^t}e_i(z;V_\rho\otimes\Delta)\overline{w(z)}e_j(z;V_\rho^\vee\otimes\Delta)\widehat{dz}\Big)\equiv 0$$ for any fixed bases $\{e_i(z;V_\rho^{(*)}\otimes\Delta)\}$ of $H^0(M,V^{(*)}\otimes\Delta).$
\end{thm}
This generalizes the standard theory of abelian theta functions and the Brill-Noether loci
to non-abelian setting. For example, when $r=1$, on the $i$-th Brill-Noether locus, the analogue of (iv) says that
the analytic torsion may be calculated via the norm of the $i$-th partial derivatives of the standard theta functions. For details, see [F, Thm 4.9] and [ACGR]. 

Put all this together,
we have then justified the convergence in the abelian case, namely, $r=1$. As for general non-abelian cases, 
such a strong result has yet  been obtained. Fortunately, what needed is a much weak result which we recall below:

\begin{thm} ({\bf Degenerations of Analytic Torisons} [F, Thm 4.12]) Let $L$ be a line bundle of degree $d\geq g-1$ such that $h^1(V_\rho\otimes L)=n>0$ for a fixed $V_\rho\in \mathcal M_{M,r}(0)$. Then within a neighborhood
$U$ of $V_\rho$ in $\mathcal M_{M,r}(0)$, the analytic torsions $T(V_{\rho(s)}\otimes L)$, which is positive whenever $h^1(V_{\rho(s)}\otimes L)=0$, vanishes to order $2n$ at $V_{\rho}=V_{\rho(0)}$. In particular, near $s=0$,
$$T(V_{\rho(s)}\otimes L)=4^nT(V_\rho\otimes L)\det[\overline{{^t}C(s)}C(s)]+O(\|s\|^{2n+1})
$$ where for any orthonormal basis $\{e_i\}, \{e_j^*\}$ for $H^0(M,V_\rho\otimes L)$, $H^0(M,K_M\otimes (V_\rho\otimes L)^\vee)$ respectively:
$$C_{ij}(s):=\int_M{^te}_i(z;V_{\rho}\otimes L)\overline w e_j^*(z;K_M\otimes (V_\rho\otimes L)^\vee)\widehat {dz}.$$ And for any $x_1,\dots, x_{d+1-g}$,
$$\begin{aligned}&\det[\overline{{^t}C(s)}C(s)]\det[B(x_i,\overline{x_j};V_{\rho(s)}\otimes L)]\\
=&\Biggl|\det\begin{pmatrix}{^te}_1(x_1:V_\rho\otimes L)&\cdots&{^te}_1(x_1:V_\rho\otimes L)&C_{11}(s)&\cdots&C_{1n}(s)\\
\cdots&&&&\\
{^te}_p(x_1:V_\rho\otimes L)&\cdots&{^te}_p(x_1:V_\rho\otimes L)&C_{p1}(s)&\cdots&C_{pn}(s)\end{pmatrix}\Biggr|^2\\
&+O(\|s\|^{2n+1})\end{aligned}$$ as $\rho(s)\to \rho$ along any smooth curve transverse at $\rho$ to the subvariety $V_1$ of all $V_\rho\in \mathcal M_{M,r}(0)$ with $h^1(M,V_\rho\otimes L)>0$.
Here $p=n+r(d+1-g)$.
\end{thm}

Thus by the fact that the Bergman kernel admits logarithmic degeneration, we complete
the proof of the convergence of the regularized integrations appeared in the proof of our zeta functions.

\subsubsection{Asymptotics of Analytic Torsions}

To establish the convergence of the infinite sum appeared in the definition
of our zeta functions on $m$, we only need to understand
$$III(s)=\sum_{m>2g-2}\int_{\mathcal M_{M,r}(rm)}\Big(e^{\tau(M,V)}-1\Big)(e^{-s})^{\chi(M,V)}d\mu,$$
by the discussion in \S1.3. Note that for stable bundle of rank $r$ and degree $mr$ with $m>2g-2$, $h^0(M,V)=\chi(M,V)=r[m-(g-1)]$ is a constant. So $T(V)$ is a constant function on $\mathcal M_{M,r}(mr)$. Thus using the natural isomorphism
$$\mathcal M_{M,r}(0)\simeq \mathcal M_{M,r}(mr),\qquad V\mapsto V\otimes A^{\otimes m},$$
we have
$$III(s)=\sum_{m>g-1}\int_{V\in \mathcal M_{M,r}(0)}\Big(e^{\tau(M,V\otimes A^{\otimes (m+g-1)})}-1\Big)(e^{-s})^{rm}d\mu.$$
As such, then the convergence is guaranteed by the following results of Faltings, Miyaoka, Bismut-Vasserot:

\begin{thm} (See e.g., [BV, Thm 8])
As $m\to\infty$,
$$\frac{d}{ds}\zeta^0_\lambda(s,V\otimes A^{\otimes (m+1-1)})\Big|_{\lambda=s=0}=
O(m^r\log m).$$
\end{thm}

Thus, in essence, we are dealing with the infinite summation
$$\sum_{m>g-1}\int_{M_{m,r}(0)}\Big(e^{e^{-m^r\log m}}-1\Big)
(e^{-s})^{rm}d\mu$$ which is clearly convergent. This then completes the proof of Thm 1 stated in \S1.2.

\section{General Uniformity of Zetas}
\subsection{Number Fields: $SL_n$}
\subsubsection{Siegel's Volume Formula}
For special linear group $SL_n$ defined over $\mathbb Q$, there are 3 naturally associated groups, namely, the real Lie group $SL_n(\mathbb R)$,
 its maximal compact subgroup $SO_n(\mathbb R)$ and the full modular group $SL_n(\mathbb Z)$. It is well known that the double quotient space $SL_n(\mathbb Z)\backslash SL_n(\mathbb R)/O_n(\mathbb R)$ may be interpreted as the space of isometric classes of rank $n$ lattices of volume one 
 in the Euclidean space $\mathbb R^n$. Indeed, the metrics on $\mathbb R^n$ are parametrized by matrices $A\cdot A^t$ with $A\in GL_n(\mathbb R)$, and up to $O_n(\mathbb R)$-equivalence, the metric is uniquely determined by $A$. As such, then the lattice structures are finally determined modulo the automorphism group $SL_n(\mathbb Z)$ of $\mathbb Z^n$. 
 
Denote by $\mathbb M_{\mathbb Q,n}[1]$ the moduli space of all full rank lattices in $\mathbb R^n$ of volume one. The above discussion exposes
the following
 
 \noindent
 \begin{fact} ({\bf Arithmetic versus Geometry}) There is a natural one-to-one correspondence
 $$SL_n(\mathbb Z)\backslash SL_n(\mathbb R)/SO_n(\mathbb R)\,\simeq\,
 \mathbb M_{\mathbb Q,n}[1].$$
\end{fact}
Associated to the natural measure on $ SL_n(\mathbb R)$, we may ask what is the corresponding volume of the above space.
Surprisingly, while the space $SL_n(\mathbb Z)\backslash SL_n(\mathbb R)/SO_n(\mathbb R)$, namely,
$\mathbb M_{\mathbb Q,n}[1]$, is highly non-abelian, or better, non-commutative, 
according to Siegel, its volume can be expressed in terms of the special values of Riemann zeta function, which is abelian in nature.
\vskip 0.30cm
\begin{fact} (Siegel) ({\bf Volume of Fundamental Domain})
$$m_{\mathbb Q,n}:=\mathrm{Vol}\Big(\mathbb M_{\mathbb Q,n}[1]\Big)=\widehat\zeta_{\mathbb Q}(1)
\widehat\zeta_{\mathbb Q}(2)\cdots
\widehat\zeta_{\mathbb Q}(n)$$
where $\widehat\zeta_{\mathbb Q}(s)$ denotes the complete Riemann zeta function and $$\widehat\zeta_{\mathbb Q}(1)
:=\mathrm{Res}_{s=1}\widehat\zeta_{\mathbb Q}(s).$$
\end{fact}
\subsubsection{Stability}
Among all lattices, motivated by Mumford's fundamental work in algebraic geometry,
we independently introduced the semi-stable lattices in our studies of non-abelian zeta functions. By definition, a lattice $\Lambda$ is called {\it semi-stable} if for all sub-lattices $\Lambda_1$ of $\Lambda$
$$\mathrm{Vol}(\Lambda_1)^{\mathrm{rank}\Lambda}\geq 
\mathrm{Vol}(\Lambda)^{\mathrm{rank}\Lambda_1}.$$
Denote by $\mathbb M_{\mathbb Q,n}^{\mathrm{ss}}[1]$ the moduli space of rank $n$ semi-stable lattices of volume 1. 
One checks that $\mathbb M_{\mathbb Q,n}^{\mathrm{ss}}[1]$ is a closed compact subset of $\mathbb M_{\mathbb Q,n}[1]$. With the induced metric, define
$$m_{\mathbb Q,n}^{\mathrm{ss}}:=\mathrm{Vol}\Big(\mathbb M_{\mathbb Q,n}^{\mathrm{ss}}[1]\Big).$$
A natural question is what is the volume $u_{\mathbb Q,n}$.

\subsubsection{High rank non-abelian zeta functions}

Similarly, denote by $\mathbb M_{\mathbb Q,n}^{\mathrm{ss}}$ the moduli space of rank $n$ semi-stable lattices  and by $\mathbb M_{\mathbb Q,n}^{\mathrm{ss}}[T]$ its volume $T$ part. 
Then we have a natural decomposition
$$\mathbb M_{\mathbb Q,n}^{\mathrm{ss}}=\cup_{T\in\mathbb R_{>0}}\mathbb M_{\mathbb Q,n}^{\mathrm{ss}}[T].$$
Easily one checks that there is a natural isomorphism
$$\mathbb M_{\mathbb Q,n}^{\mathrm{ss}}[T]\simeq\mathbb M_{\mathbb Q,n}^{\mathrm{ss}}[T']\qquad\forall\ T,\,T'\in\mathbb R_{>0}.\eqno(*)$$ Using the above measure on $\mathbb M_{\mathbb Q,n}^{\mathrm{ss}}[T]$ and the invariant Haar measure $\frac{dT}{T}$ on $\mathbb R_{>0}$, we obtain a natural measure $d\mu$ on $\mathbb M_{\mathbb Q,n}^{\mathrm{ss}}$.

Moreover, there is a genuine cohomology theory $h^i(F,\Lambda), i=0,1$ for lattices $\Lambda$ over number fields $F$ for which the arithmetic analogue of the duality, the Riemann-Roch theorem, the vanishing theorem holds. 
For details, please refer to [W]. In the case of $F=\mathbb Q$, $$h^0(\mathbb Q,\Lambda)=\log\Big(\sum_{{\bf x}\in\Lambda}e^{-\pi\|{\bf x}\|^2}\Big)$$ which was introduced earlier in [GS]. Denote by $d(\Lambda)$ the Arakelov degree of $\Lambda$, which over $\mathbb Q$ is simply $-\log\mathrm{Vol}(\Lambda)$.
Following [W], define the associated {\it rank $n$ non-abelian zeta function}
$\widehat\zeta_{\mathbb Q,n}(s)$ by
$$\widehat\zeta_{\mathbb Q,n}(s):=\int_{ \mathbb M_{\mathbb Q,n}^{\mathrm{ss}}}\Big(e^{h^0(\mathbb Q,\Lambda)}-1\Big)\cdot(e^{-s})^{d(\Lambda)}\,d\mu,\qquad\mathrm{Re}(s)>1.$$
Then using the basic property of the above cohomology theory for $h^i$'s, namely the duality, the RR and the vanishing theorem,
tautologically, we have the following
\begin{fact} (Weng) (0) ({\bf Relation with Abelian Zeta})  $$\widehat\zeta_{\mathbb Q,1}(s)=\widehat\zeta_{\mathbb Q}(s);$$
(i)  ({\bf Meromorphic Extension}) $\widehat\zeta_{\mathbb Q,n}(s)$ is a well-defined holomorphic function in $\mathrm{Re}(s)>1$, and admits a unique meromorphic extension to the whole $s$-plane;

\noindent
(ii) ({\bf Functional Equation})
$$\widehat\zeta_{\mathbb Q,n}(1-s)=\widehat\zeta_{\mathbb Q,n}(s);$$

\noindent
(iii) ({\bf Singularities)}
$\widehat\zeta_{\mathbb Q,n}(s)$ has only two singularities, all simple poles, at $s=0,\,1$. Moreover
$$\boxed{\mathrm{Res}_{s=1}\widehat\zeta_{\mathbb Q,n}(s)=m_{\mathbb Q,n}^{\mathrm{ss}}:=\mathrm{Vol}\Big(\mathbb M_{\mathbb Q,n}^{\mathrm{ss}}[1]\Big)}$$
\end{fact}

In particular we see that $m_{\mathbb Q,n}^{\mathrm{ss}}$ is naturally related to the special value of the non-abelian zeta function $\widehat\zeta_{\mathbb Q,n}(s)$.
\subsubsection{Parabolic Reduction: Analytic Theory}
The high rank zeta functions are closely related with Eisenstein series. In fact, we have
\begin{fact} (Weng) (i) ({\bf High Rank Zeta and Eisenstein Series})
$$ \widehat\zeta_{\mathbb Q,n}(s)=\int_{\mathbb M_{\mathbb Q,n}^{\mathrm{ss}}[1]}\widehat E(\Lambda,s)\,d\mu
=\int_{\Big(SL_n(\mathbb Z)\backslash SL_n(\mathbb R)/O_n(\mathbb R)\Big)^{\mathrm{ss}}}\widehat E^{SL_n/P_{n-1,1}}({\bf 1},g;s).$$ Here 
$\widehat E(\Lambda,s)$ denotes the complete Eisenstein series associated to the lattice $\Lambda$, $\widehat E^{SL_n/P}({\bf 1},g;*)$
denote the relative (complete) Eisenstein series on $SL_n(\mathbb R)$ induced from the constant function ${\bf 1}$ on the Levi
 factor of the maximal parabolic subgroup $P$ and $P_{n-1,1}$ denotes the standard parabolic subgroup of $SL_n$ corresponding
 to the partition $n=(n-1)+1$, and $\Big(SL_n(\mathbb Z)\backslash SL_n(\mathbb R)/O_n(\mathbb R)\Big)^{\mathrm{ss}}$ the part corresponding to the semi-stable lattices via Fact 1, which for our convenience  will also be viewed as a subset of $SL_n(\mathbb R)$;

\noindent
(ii) ({\bf Analytic Truncation versus Arithmetic Truncation})
 $$\Lambda^{0}{\bf 1}=\chi_{\Big(SL_n(\mathbb Z)\backslash SL_n(\mathbb R)/O_n(\mathbb R)\Big)^{\mathrm{ss}}}$$
Namely, Arthur's truncation of the constant function ${\bf 1}$ is simply  the characteristic function of the subset $\Big(SL_n(\mathbb Z)\backslash SL_n(\mathbb R)/O_n(\mathbb R)\Big)^{\mathrm{ss}}$ consisting of semi-stable points.
  \end{fact}
This is a number theoretic analogue of a result of Laffourge ([L]) on the relation
between analytic truncation and arithmetic truncation for function fields.

Consequently, $$ \widehat\zeta_{\mathbb Q,n}(s)
=\int_{SL_n(\mathbb Z)\backslash SL_n(\mathbb R)/O_n(\mathbb R)}\Lambda^{\bf 0}\widehat E^{SL_n/P_{n-1,1}}({\bf 1},g;s).$$
Here $\Lambda^{\bf 0}\widehat E^{SL_n/P_{n-1,1}}({\bf 1},g;s)$ denotes the Arthur's truncation of the Eisenstein series $\widehat E^{SL_n/P_{n-1,1}}({\bf 1},g;s)$.
On the other hand, by Langlands' theory of Eisenstein series, we know that
$$\widehat E^{SL_n/P_{n-1,1}}({\bf 1},g;s)=\mathrm{Res}_{\langle \lambda-\rho,\alpha_i^\vee\rangle=0, i=1,2,\dots,n-2}
\widehat E^{SL_n/P_{1,\dots,1}}({\bf 1},g;\lambda)$$ where $\alpha_i=\alpha_i-\alpha_{i+1}$ denotes the simple roots of the root system $A_{n-1}$ associated to $SL_n$, and $\rho=\frac{1}{2}\sum_{\alpha>0}\alpha$ the Weyl vector.

With this, now notice that the moduli space $\mathbb M_{\mathbb Q,n}[1]$ is compact, and that on the Levi of the Borel subgroup, ${\bf 1}$ is cuspidal. So we can evaluate the {\it Eisenstein period}
$$\int_{SL_n(\mathbb Z)\backslash SL_n(\mathbb R)/O_n(\mathbb R)}\Lambda^{\bf 0}\widehat E^{SL_n/P_{1,\dots,1}}({\bf 1},g;\lambda).$$ This then gives a very precise expression of non-abelian zeta function $\widehat\zeta_{\mathbb Q,n}(s)$ as a combination of terms consisting of products of rational functions coming from the symmetry depending only on the root system, and  abelian zeta functions.
For details, please see [W]. As a direct consequence, we have the following
\begin{fact} (Weng) ({\bf Parabolic Reduction, Stability \& the Volumes})
$$m_{\mathbb Q,n}^{\mathrm{ss}}=\sum_{k\geq 1}(-1)^{k-1}\sum_{n_1+\dots+n_k=n,n_i>0}\frac{1}{\prod_{j=1}^{k-1}(n_j+n_{j+1})}\cdot \prod_{j=1}^km_{\mathbb Q,n_j}.$$
\end{fact}
Geometrically, this means that the semi-stable part can be obtained from the fundamental domain associated to $SL_n$ by deleting the tubular neighborhoods of cusps corresponding to parabolic subgroups which parametrize the same type of canonical flags  of unstable lattices,
and whose volumes, up to the lattice extensions, are completely determined by that associated to the simple factors of related Levi factors. Undoubtedly, this parabolic reduction is also the \lq heart' of the theory of the truncations, both, analytic and arithmetic.
\subsubsection{Parabolic Reduction: Geometric Theory}
During our Sept, 2012's stay at IHES, Kontsevich introduced us their beautiful formula relating $m_{\mathbb Q,n}$'s and $m_{\mathbb Q,n}^{\mathrm{ss}}$'s. This basic relation is obtained within their lecture notes on the wall-crossing ([KS]), with suitable
renormalizations, their result can be stated as follows:

\begin{fact} (Kontsevich-Soibelman) ({\bf Parabolic Reduction, Stability \& the Volume})
$$\frac{1}{n}\cdot m_{\mathbb Q,n}=\sum_{k\geq 1}\sum_{n_1+\dots+n_k=n,n_i>0}c_{n_1,n_2,\dots,n_k}\cdot \prod_{j=1}^km_{\mathbb Q,n_j}^{\mathrm{ss}}.$$ Here $\displaystyle{c_{n_1,n_2,\dots,n_k}:=\frac{1}{n_1(n_1+n_2)\cdots(n_1+n_2+\cdots+n_k)\cdots(n_{k-1}+n_k)n_k}}.$
\end{fact} 
Indeed, the essence of this is the existence of the so-called canonical filtration, namely, the Harder-Narasimhan filtration of a lattice: {\it For a rank $n$ lattice $\Lambda$, there exists a unique filtration of sub-lattices
$$0=\Lambda_0\subset\Lambda_1\subset\Lambda_2\subset\cdots\subset\Lambda_k=\lambda$$ such that}

(i) {\it $G_i(\Lambda):=\Lambda_i/\Lambda_{i+1}$ is semi-stable}; and

(ii) $\displaystyle{\mathrm{Vol}\Big(G_i(\Lambda)\Big)^{\mathrm{rank}\big(G_{i+1}(\Lambda)\big)}> \mathrm{Vol}\Big(G_{i+1}(\Lambda)\Big)^{\mathrm{rank}\big(G_{i}(\Lambda)\big)}}.$

\subsection{Function Fields$/\mathbb F_q$: $SL_n$}
The interactions between studies of number fields and that of function fields (over finite fields)
have been proven to be very fruitful, based on formal analogues between these two types of fields, despite the fact that many working mathematicians, not without their own reasons, believe otherwise. The works  presented
here are yet another group of beautiful examples.
\subsubsection{Weil's Formula: Tamagawa Numbers}
Motivated by Siegel's volume formula above, which originally was done in the theory of quadratic forms, Weil reinterpreted it in terms his famous Tamagawa number one conjecture ([Weil]). For $SL_n$, this goes as follows.

Let $X$ be an irreducible, reduced, regular projective curve defined over $V$. 
Denotes its function field by $F$ and its ring of adeles by $\mathbb A$. Fix a vector bundle $E_0$ of rank $n$ on $X$ with determinant $\lambda$. 

Consider the group $SL_n(\mathbb A)$ with $\mathbb K(E_0)$ the maximal compact subgroup associated to $E$. Then there is a natural morphism $\pi$ from the quotient space $SL_n(F)\backslash SL_n(\mathbb A)/\mathbb K(E_0)$ to the stack $\mathbb M_{X,n}(\lambda)$ of rank $n$ bundles with fixed determinant $\lambda$. 

\begin{fact}  ([HN], [DR]) The natural morphism $$\pi:SL_n(F)\backslash SL_n(\mathbb A)/\mathbb K(E_0)\to\mathbb M_{X,n}(\lambda)$$
is surjective with the fiber $\pi^{-1}(E_0^g)$ at the vector bundle
$E_0^g$ associated to $g\in SL_n(\mathbb A)$ consisting of 
$\#\Big(\mathbb F_q^*/\det\mathrm{Aut}(E_0^g)\Big)$. Here $\det\mathrm{Aut}(E_0^g)$ denotes the image of $\det\mathrm{Aut}(E_0^g)$ in $\mathbb F_q^*$ under the determinant mapping.
\end{fact}

Denote by $\mathbb M_{X,n}(d)$ the moduli stack of rank $n$ bundle of degree $d$ on $X$, and introduce the total mass for rank $n$ and degree $d$ bundles on $X$ by $$m_{X,n}(d):=\sum_{E\in \mathbb M_{X,n}(d)}\frac{1}{\#\mathrm{Aut}(E)}.$$ Denote by
$$
\widehat\zeta_X(s):=\zeta_X(s)\cdot (q^s)^{g-1}$$  the complete Artin zeta function associated to $X$, and $$\widehat\zeta_X(1):=\mathrm{Res}_{s=1}\widehat\zeta_X(s)\cdot \log q.$$
Then Weil's result that the Tamagara number of the quotient space $SL_n(F)\backslash SL_n(\mathbb A)$ equals one is equivalent to the following
\begin{fact} (Weil) ({\bf Tamagawa Number})
$$ m_{X,n}(d)=m_{X,n}:= \widehat\zeta_X(1)
\widehat\zeta_X(2)\cdots \widehat\zeta_X(n).$$ 
\end{fact}

\subsubsection{Non-Abelian Zeta Functions for $X$}
Denote by  $\mathbb M_{X,n}^{\mathrm{ss}}(d)$ the moduli stack of rank $n$ semi-stable bundle of degree $d$ on $X$. Then define the {\it pure non-abelian zeta function of rank $n$ for $X$} by
$$\widehat\zeta_{X,n}(s):=\sum_{k\in\mathbb Z}\sum_{E\in \mathbb M_{X,n}^{\mathrm{ss}}(kd)}\frac{q^{h^0(X,E)}-1}{\#\mathrm{Aut}(E)}\cdot (q^{-s})^{\chi(X,E)}.$$ 
Write $$\widehat\zeta_{X,n}(s)=\zeta_{X,n}(s)\cdot (q^s)^{n(g-1)},\qquad Z_{X,n}(t):=\zeta_{X,n}(s)\ \mathrm{with}\ t=q^{-s}$$

Introduce the partial mass of  semi-stable bundles by
$$\alpha_{X,n}(d):=\sum_{E\in \mathbb M_{X,n}^{\mathrm{ss}}(d)}\frac{q^{h^0(X,E)}-1}{\#\mathrm{Aut}(E)},\qquad
\beta_{X,n}(d):=\sum_{E\in \mathbb M_{X,n}^{\mathrm{ss}}(d)}\frac{1}{\#\mathrm{Aut}(E)}.$$
Then tautologically, 
$$\begin{aligned} Z_{X,n}&(t)
=\sum_{m=0}^{(g-1)-1}\alpha_{X,n}(mn)
\cdot \Big(T^{m}+ Q^{(g-1)-m}\cdot T^{2(g-1)-m}\Big)\\
&+\alpha_{X,n}\big(n(g-1)\big)\cdot T^{g-1}+\frac{\beta_{X,r}(0)T^{2g-1}}{(1-QT)(1-T)}\cdot\Big[(Q^g-1)-(Q^g-Q)T\Big].
\end{aligned}$$ where $T:=t^n$ and $Q:=q^n$. This exposes the following

\begin{fact} (Weng) (i) ({\bf Relation with Artin Zetas})
$\zeta_{X,1}(s)=\zeta_X(s)$, the Artin zeta function for $X/\mathbb F_q$;

\noindent
(ii)  ({\bf Rationality}) There exists a degree $2g$ polynomial $P_{X,r}(T)\in\mathbb Q[T]$ of $T$
such that
$$Z_{X,r}(t)=\frac{P_{X,r}(T)}{(1-T)(1-QT)}\quad\mathrm{with}\quad T=t^r,\ Q=q^r;$$
\noindent
(iii)  ({\bf Functional Equation}) $$\widehat \zeta_{X,n}(1-s)=\widehat \zeta_{X,n}(s);$$

\noindent
(iv) ({\bf Residues}) $$\widehat \zeta_{X,n}(1):=\mathrm{Res}_{s=1}\widehat \zeta_{X,n}(s)\cdot {\log Q}=\beta_{X,n}(0)\Big(=m_{X,n}^{\mathrm{ss}}(0)\Big).$$
\end{fact}

\subsubsection{Parabolic Reduction, Stability and the Mass: Geometric Theory}

To go further, make a normalization by introducing
 $$\widetilde m_{X,n}^{\mathrm{ss}}(d)=\frac{1}{q^{\frac{n(n-1)}{2}(g-1)}}\cdot m_{X,n}^{\mathrm{ss}}(d).$$
Then the parabolic reduction via the Harder-Narasimhan filtration leads to the following relation involving infinite sums:

\begin{fact} ([HN], [DR]) ({\bf Parabolic Reduction})
$$ m_{X,n}(d)=\sum_{k\geq 1}\sum_{n_1+\dots+n_k=n,n_i>0}\sum_{\substack{\frac{d_1}{n_1}>\dots>\frac{d_k}{n_k}\\
d_1+\dots+d_k=d}}q^{-\sum_{i<j}(d_in_j-d_jn_i)}\prod_{j=1}^k\widetilde m_{X,n_j}^{\mathrm{ss}}(d_j).$$
\end{fact}

\subsubsection{Parabolic Reduction, Stability and the Mass: Combinatorial Aspect}
With the above result, Zagier proved the following fundamental result, hidden in his paper on Verlinder formula:

\begin{fact}([Z], see also [WZ2]) ({\bf Parabolic Reduction, Stability \& the Mass})
$$\widetilde m^{\mathrm {ss}}_{X,n}(d)=\sum_{k\geq 1}(-1)^{k-1}\sum_{\substack{n_1+\dots+n_k=n\\ n_i>0, i=1,\dots,k}}\prod_{j=1}^{k-1}\frac{q^{(n_j+n_{j+1})\cdot\{(n_1+\dots+n_j)\cdot\frac{d}{n}\}}}{q^{(n_j+n_{j+1})}-1}\cdot\prod_{j=1}^k m_{X,n_j}.$$
\end{fact} 

This formula should be compared with with our Fact 5 for number fields. The structure are very much similar: in fact if we let $q\to 1$, then we would get the number theoretic identity there.

\subsubsection{Parabolic Reduction, Stability and the Mass: New Formula}
The above relation of Harder-Narasimhan, Ramanan-Desale and Zagier for function fields correspond to our own formula listed as Fact 5. So naturally, what should be the one appeared in the theory of
parabolic reduction, stability and the volumes obtained by Kontsevich-Soibelman ([KS]).

\begin{fact}  ([WZ2]) For an ordered partition $n=n_1+n_2+\dots+n_k$, fix $\delta_i\in \{0,,\dots, n_1-1\}$, then fix $v_i\in [0,1)\cap\mathbb Q$ satisfying $$v_i\equiv\frac{\delta_i}{n_i}-\frac{\delta_{i+1}}{n_{i+1}}\ (\mathrm{mod}\, 1).$$ Fix an $n$-th primitive root of unity $\zeta_n$.
Also set 
$N_i=n_1+n_2+\dots+n_i$ and $N_i'=n-N_i$ for $i=1,2,\dots n$.

\noindent
(i) ({\bf On Average})
$$\begin{aligned}n\cdot m_{X,n}=&\sum_{k\geq 1}\sum_{\substack{n_1+n_2+\dots+n_k=n\\ n_i>0,\, i=1,2,\dots,k}}
\sum_{\substack{\delta_i\in\{0,1,\dots,n_i-1\}\\j=1,2,\dots,k}}
\prod_{i=1}^{k-1}\frac{q^{v_iN_iN_i'}}{q^{N_iN_i'}-1}\cdot 
\prod_{j=1}^k\widetilde m^{\mathrm{ss}}_{X,n_j}(\delta_j).
\end{aligned}$$

\noindent
(ii) ({\bf Individuality}) For all $d=0,1,\dots,n-1$,
$$\begin{aligned}m_{X,n}=&\sum_{k\geq 1}\sum_{\substack{n_1+n_2+\dots+n_k=n\\ n_i>0,\, i=1,2,\dots,k}}
\frac{ 1}{n}\\
&\times
\sum_{\substack{\delta_i\in\{0,1,\dots,n_i-1\}\\j=1,2,\dots,k}}
\Bigg(\sum_{\zeta_n:\,\zeta_n^n=1}\zeta_n^{n-d}\cdot \prod_{h=1}^{k-1}\frac{\zeta_n^{v_hN_h}q^{v_hN_hN_h'}}{
\zeta_n^{N_h}q^{N_hN_h'}-1}\Bigg)\prod_{j=1}^k\widetilde m^{\mathrm{ss}}_{X,n_j}(\delta_j).\end{aligned}$$
\end{fact}
The first is obtained by taking average on $d$ from the relation of Fact 11, while the second is obtained directly from that of Fact 11. With geometric
picture in mind, formula (ii) should be further polished, so as to get everything done according to the real world structure. This may prove to be a bit complicated due to the fact that usually 
$$\mathbb M_{X,n}(d)\not\simeq \mathbb M_{X,n}(d'),\qquad \forall d,\,d'\in\{0,1,\dots,n-1\},\ d\not= d'.$$
This is very different from the cases for  number fields, where we always have the isomorphism $(*)$ between different levels.

We remind the reader that while all relations in function field case are obtained using geometric methods, our basic relation for number fields are obtained analytically using Eisenstein series.

\subsection{Parabolic Reduction, Stability and the Mass: General Reductive Groups}
\subsubsection{Parabolic Reduction Conjecture}
Motivated by the above discussion, more generally, for  a split reductive group $G$  defined over a number field $F$, $B$ a fixed Borel, ...,
denote by $G(\mathbb A)^{\mathrm{ss}}$ the adelic elements of $G$ corresponding to semi-stable principle $G$-lattices ([G]). Write $\mathbb K_G$ for the associated maximal compact subgroup.
 Also for a standard parabolic subgroup $P$, write its  Levi decomposition as  $P=UM$ with $U$ the unipotent radical and $M$ its Levi factor. Denote the corresponding simple decomposition of $M$ as
$\prod_iM_i$ with $M_i$'s the simple factors of $M$.
Introduce invariants $$m_{F;P}:=\prod_i\mathrm{Vol}\Big(\mathbb K_{M_i} \big\backslash M_i^1(\mathbb A)\big/M_i(F)Z_{M_i^1(\mathbb A)}\Big)$$
and $$m_{F;P}^{\mathrm{ss}}:=\prod_i\mathrm{Vol}\Big(\mathbb K_{M_i} \big\backslash M_i^1(\mathbb A)^{\mathrm{ss}}\big/M_i(F)Z_{M_i^1(\mathbb A)}\Big).$$ In parallel, we have similar constructions for function fields $F=\mathbb F_q(X)$. 

Denote by $$n_i:=\#\{\alpha>0:\langle\rho,\alpha^\vee\rangle=i\}-
\#\{\alpha>0:\langle\rho,\alpha^\vee\rangle=i+1\}$$ and by $v_G$
the volume of $\{\sum_{\alpha\in\Delta}:a_\alpha\alpha^\vee:a_\alpha\in [0,1)\}$. 
\begin{fact} (Langlands) ({\bf Volume of Fundamental Domain}) For the field of rationals,
$$\mathrm{Vol}\Big(\mathbb K_G \big\backslash G^1(\mathbb A)\big/G(\mathbb Q)Z_{G^1(\mathbb A)}\Big)=v_G\cdot
\prod_{i\geq 1}\widehat\zeta(i)^{-n_i}.$$
\end{fact}

Based on all this, then we have the following

\begin{conj} ({\bf Parabolic Reduction}) Let $G/F$ be a split reductive group with $B/F$ a fixed Borel. Then, for each standard parabolic subgroup $P$ of $G$, there exist constants 
$c_P\in\mathbb Q,\ e_P\in\mathbb Q_{>0}$, independent of $F$,  such that
 $$m_{F;G}=\sum_{P}\,c_P\cdot m^{\mathrm{ss}}_{F;P},\qquad
 m^{\mathrm{ss}}_{F;G}=\sum_{P}\,\mathrm{sgn}(P)\cdot e_P\cdot m_{F;P},$$ where $P$ runs over all standard parabolic subgroups of $G$, and $\mathrm{sgn}(P)$ denotes the sign of $P$.  
 \end{conj}
 
\noindent
{\it Remark.} Calculations in [Ad] for lower rank groups indicate that, for number fields, $\frac{1}{c_P}\in\mathbb Z_{>0}$. It would be very interesting to find a close formula for them.
 \subsubsection{Precise Formulation: Number fields}
The exact values of $e_P$'s can be written out in terms of the associated root system. 
Indeed, 
if $$W_0:=\Big\{w\in W:\{\alpha\in\Delta:w\alpha\in \Delta\cup \Phi^- \}=\Delta\Big\},$$
then there is a natural one-to-one correspondence between $W_0$ and the set of subsets of $\Delta$, and hence to the set of standard parabolic subgroups of $G$. Thus we will write
$$W_0:=\Big\{w_P:\ P\ \mathrm{standard\ parabolic\ subgroup}\Big\},$$
and, for $w=w_P\in W_0$, write  $J_P\subset\Delta$ the corresponding subset. 

\begin{conj}  Let $G$ be a connected split reductive group 
with $P$ its maximal parabolic subgroup defined over a number field $F$.

\noindent
(i)  ({\bf Relation to Zetas with Symmetries})  
$$m^{\mathrm{ss}}_{F;G}=\mathrm{Res}_{s=-c_P}\widehat\zeta_F^{(G,P)}(s)=\mathrm{Res}_{\lambda=\rho}\omega_F^{G}(\lambda);$$

\noindent
(ii)  ({\bf Parabolic Reduction, Stability \& the Volumes})
$$ m^{\mathrm{ss}}_{F;G}=\sum_{P}\frac{(-1)^{\mathrm{rank}(P)}}{\prod_{\alpha\in\Delta\backslash w_PJ_P}(1-\langle w_P\rho,\alpha^\vee\rangle)}\cdot m_{F;P},$$ where $P$ runs over all standard parabolic subgroups of $G$. 
\end{conj}

\subsubsection{Precise Formulation: Function fields}
In parallel, we have the following

\begin{conj} For a split reductive group $G$ with $P$ its maximal parabolic subgroup, both
defined over the function field $\mathbb F_q(X)$ of an
irreducible reduced regular projective curve $X$ over $\mathbb F_q$, 

\noindent
(i) ({\bf Relation to Zetas with Symmetries}) 
$$\log q\cdot  m^{\mathrm{ss}}_{F;G}(0)=\mathrm{Res}_{s=-c_P}\widehat\zeta_X^{(G,P)}(s)=\mathrm{Res}_{\lambda=\rho}\omega_X^{G}(\lambda);$$

\noindent
(ii) ({\bf Parabolic Reduction, Stability \& the Mass})
$$m^{\mathrm{ss}}_{F;G}(0)=\sum_{P}\frac{(-1)^{\mathrm{rank}(P)}}{\prod_{\alpha\in\Delta\backslash w_PJ_P}(q^{-\langle w_P\rho,\alpha^\vee\rangle+1}-1)}\cdot m_{F;P},$$ where $P$ runs over all standard parabolic subgroups of $G$. 
 \end{conj}
More generally, also to make this compatible with the paper of Laumon-Rapoport ([LR]), and for the purpose of understanding what
 happens for Riemann surfaces to be discussed below, next we consider arbitrary slope.

So in the follows we will freely use the notations of [LR]. Consider then the adelic space
$$\frak K Z_{G(\mathbb A)}\backslash G(\mathbb A)^1/G(K)$$ with $\frak K$ the maximal compact subgroup of $G(\mathbb A)$, where $\mathbb A$ denotes the adelic ring of $K$. For each $[{\bf g}]$ in this quotient space with ${\bf g}\in G(\mathbb A)^1$, (up to universal constant factors independent of ${\bf g}$,)
define its mass by $$m([{\bf g}]):=
\frac{1}{|{\bf g}\frak K{\bf g}^{-1}\cap G(K)|}\buildrel\cdot\over=\frac{1}{|\Aut(\mathcal E_{\bf g})|},$$ with $\mathcal E_{\bf g}$ the associated principal $G$-bundle, and the total mass
by $$m_{K;G;X}:=\sum_{[{\bf g}]}m([{\bf g}])$$ where ${\bf g}$ runs over all quotient classes.
Similarly, for each $\nu_G'\in X_*(A_G')$, set $m_{K;G;X}^{\mathrm{ss}}(v\nu_G')$ be the partial (total) mass defined by 
running $[{\bf g}]$ over these whose associated $G$-bundles on $X$ are semi-stable of slope $\nu_G'$.

\begin{conj} (i) (Total Mass) $$m_{K;G}(\nu_G')=m_{K;G}:=\prod_{i\geq 1}\widehat\zeta_X(i)^{-n_i}.$$

\noindent
(ii) ({\bf Parabolic Reduction}: Infinite Form)
$$m_{K;G}=\sum_{P\in\mathcal P}\sum_{\substack{\nu_P'\in X_*(A_P')\\ [\nu_P']_G=\nu_G'}}\tau_P^G([\nu_P']^G)q^{2m(P,\nu_P')}m_{K;P}^{\mathrm{ss}}(\nu_P').$$ Here as above, $m_{K;P;X}^{\mathrm{ss}}(\nu_P'):=\prod_{i=1}^km_{K;M_{P,i}}^{\mathrm{ss}}({\nu^{\dag\,'}_{P,i}})$ with $
{\nu^{\dag\,'}_{P,i}}$ the induced $M_{P,i}$ component of $\nu_P'$.

\noindent
(iii) ({\bf Parabolic Reduction}: Finite Form)
$$m_{K;G}^{\mathrm{ss}}(\nu_G')=\sum_{P\in\mathcal P}(-1)^{\dim\frak a_P^G}\sum_{\alpha\in\Delta_P}\frac{\prod_{\substack{\lambda\in\Lambda_P^G\\ [\lambda]_G=\nu_G'}}q^{2\langle\rho_P,\alpha^\vee\rangle
\langle\varpi_\alpha^G(\lambda)\rangle}}{q^{2\langle \rho_P,\alpha^\vee\rangle}-1}\cdot m_{K;P}$$ where as above $m_{K;P}:=\prod_{j=1}^k m_{K;M_{P,i}}.$
\end{conj}

Clearly, (i) is the analogue of Siegel-Langlands' results on the volumes of fundamental domains associated to $G$ over number fields, (ii) and (iii) coincides with Harder-Narasimhan ([HN]),
 Desale-Ramanan ([DR]), and Zagier ([Z], see also [WZ2]), respectively when $G=SL_r$. Following [LR], with Langlands' lemma, (iii) is a direct consequence of (i) and (ii). But (ii) is the argument of existence of Harder-Narasimhan filtration. So with ([RR]), what left is essentially (i).

As a direct consequence, then we can have Thms 3.2, 3.3 and 3.4 of [LR] on Poincare series for various moduli stacks of $G$-bundles, obtained by combining the method of Harder-Narasimhan ([HN]) and Atiyah-Bott ([AB]). 

\subsection{Riemann Surfaces}

\subsubsection{Zetas with Symmetries}
By comparing the Parabolic Reduction for number fields and function fields, and note that the natural $\infty$-factor corresponding to $\frac{1}{1-q^{-s}}$ is $\Gamma_{\mathbb R}(s):=\pi^{-\frac{s}{2}}\Gamma(\frac{s}{2})$,
we propose the following definition of zetas for Riemann surfaces with symmetries.

Set $\widehat\zeta_M(s):=\widehat\zeta_{M,1}(s)$ and call it the zeta function of $M$.
Then for any pair $(G,P)$ with $G$ a reductive group and $P$ its maximal parabolic subgroup, all defined over $\mathbb C(M)$, the period of $G$ over $M$ is defined by
$$\omega_X^G(\lambda):=\sum_{w\in W}(-1)^{|G|}\prod_{\alpha\in\Delta}
\Gamma_{\mathbb R}\Big(-\langle w\lambda-\rho,\alpha^\vee\rangle\Big)
\cdot\prod_{\alpha>0,w\alpha<0}
\frac{\widehat\zeta_X(\langle\lambda,\alpha^\vee\rangle)}
{\widehat\zeta_X(\langle\lambda,\alpha^\vee\rangle+1)},$$
and the period of $(G,P)$ over $M$ as
$$\omega_M^{(G,P)}(\lambda_p):=\mathrm{Res}_{\substack{\langle \lambda-\rho,\beta_i^\vee\rangle=0\\
\{\beta_1,\dots,\beta_{|P|}\}\subset\Delta \leftrightarrow P}}
\Big(\omega_X^G(\lambda)\Big).$$
This then finally leads to the zeta functions with symmetries associated to $(G,P)$ for $M$ as $$\widehat\zeta_M^{G/P}(s):=\mathrm{Norm}\Big(\omega_M^{(G,P)}(\lambda_p)\Big)$$ where Norm stands the normalization obtained similarly as for number fields and/or function fields.

\subsubsection{General Uniformity of Zetas: Special One}

Motivated by our works on number fields and function fields, we formulate the following:

\begin{conj} ({\bf Special Uniformity}) There exist constants $c_n, a_n,b_n$ such that
$$\widehat\zeta_{M,r}(s)=c_n\cdot \widehat\zeta_M^{SL(n)/P_{n-1,1}}(a_ns+b_n).$$
\end{conj}

\subsubsection{Two Dimensional Gauge Theory}
For the volume of the moduli spaces $\mathcal M_{M,r}(0)$, Witten in his two dimensional gauge theory established the following

\begin{thm}({\bf Volume Formula}, See e.g. [Wi])
$$\mathrm{Vol}\Big(\mathcal M_{X,n}(0)\Big)=n\cdot\Big(\frac{\mathrm{Vol}(SU(n))}{(2\pi)^{n^2-1}}\Big)^{2g-2}\sum_{\rho}\frac{1}{\mathrm{dim}\rho^{2g-2}}$$ where
$\rho$ runs over all irreducible representations of $SU(n)$.
\end{thm}
In particular, we note that these volumes are independent of the Riemann surfaces. 
Such an independence is in fact embedded in the main result of conformal field
theory claiming that conformal blocks, namely, the global sections of multiples of our 
determinant line bundle $\lambda$ above,
 are essentially independent of the complex structures used: after all, they are projectively flat. Thus, particularly, their dimensions, and hence its leading term, namely our volume up to normalization, according to the Hirzebruch-Riemann-Roch, are independent of  the complex structures used.
 \subsubsection{General Uniformity of Zetas: Tamagawa Number Conjecture}
 \begin{itemize}
\item $F=\mathbb C(X)$: function field of $X$
\item $\mathbb A$: ring of adeles of $F$
\item $E_0:$ fixed bundle $E_0$ w/ determinant $\lambda$, say $E=\oplus_{i=0}^{n-1}\oplus \lambda$
\item $\mathbb K=\prod_{p\in X}\mathbb K_p$
w/ $\mathbb K_p:=SL_n(\mathbb C[[z_p]])$
\end{itemize}

\begin{fact}({\bf Adelic Interpretation of Space of Bundles})
The space of rank $n$ bundles may be identified with $$\Big(\prod_pGL_n(\mathbb C[[z_p]])\Big)\Big\backslash GL_n(\mathbb A)\Big/GL_n(F)$$
\end{fact}
 Moreover, we have the following
 
 \begin{conj} ({\bf Tamagawa Numbers})
 There exist natural measures such that

\noindent
(i) $\tau\Big(SL_n(F)\backslash SL_n(\mathbb A)\Big)=1$

\noindent
 (ii) $\mathrm{m}(\mathbb K)=\widehat\zeta_M(1)\,\widehat\zeta_M(2)\cdots \widehat\zeta_M(n)$.
\end{conj}
 \subsubsection{General Uniformity of Zetas: Parabolic Reduction}
 Here the difficulty is to count $\Aut(\mathcal E)$ for a vector bundle $\mathcal E$ on $M$: even for stable bundles, it gives $\mathbb C^*=GL_1(\mathbb C)$. Any natural count for $GL_n(\mathbb C)$, or more generally, the space $\mathbb C^n$ and the torus $(\mathbb C^*)^m$.
 to make the discussion go forward, assume that we do have a natural way to count all of them (in a compatible way), then we would get
 a nice definition of the mass for  each bundle $\mathbb E$ and hence also $m(SL_n;M)$, the corresponding
 $m^{\mathrm{ss}}(SL_n;M)$, and moreover $m(SL_n;P)$, $m^{\mathrm{ss}}(SL_n;P)$
 for each parabolic subgroup $P$. By a parallel way of thinking, we then would arrive as the following
 
 \begin{conj} ({\bf General Uniformity of Zetas}: Riemann Surfaces) For a connected compact Riemann surface $M$,
 
 \noindent
 (i) (Total Mass) $$m_{M;G}(\nu_G')=m_{M;G}:=\prod_{i\geq 1}\widehat\zeta_M(i)^{-n_i}.$$

\noindent
(ii) ({\bf Parabolic Reduction}: Finite Form)
$$m_{M;G}^{\mathrm{ss}}(\nu_G')=\sum_{P\in\mathcal P}(-1)^{\dim\frak a_P^G}\sum_{\alpha\in\Delta_P}
\Gamma_{\mathbb R}(\langle \rho_P,\alpha^\vee\rangle-1)\cdot m_{M;P}.$$ \end{conj}

\subsubsection{A Rational Function}

Another way to approach is to use $h^0$ directly to create non-abelian zeta functions for Riemann surfaces. This goes as follows:

Let  $X$ be a compact Riemann surface of genus $g$. Denote by
$\mathcal M^{\mathrm{ss}}_{X,n}(d)$ the moduli space of semi-stable bundles of rank $n$ and degree $d$ on $X$
and $W_{X,n}^{\geq i}:=\{V\in \mathcal M^{\mathrm{ss}}_{X,n}(d): h^0(X,V)\geq i\}$ the Brill-Noether filtration
with $W_{X,n}^{=i}:=W_{X,n}^{\geq i}\backslash W_{X,n}^{\geq i+1}$. By the above discussion, we get
$d\mu:=d\mu_{X,n}(d)$ naturally induced volume forms on $W_{X,n}^{=i}(d)$.

Fix $\pi\in\mathbb R$, introduce the
 {\it regularized integral}
$$\begin{aligned}&\int_{\mathcal M_{X,n}(d)}^{\#}\Big(\pi^{h^0(X,V)}-1\Big)\Big(\pi^{-s}\Big)^{\chi(X,V)}d\mu_{X,n}(d)\\
:=&\sum_{i\geq 0}\int_{W_{X,n}^{=i}(d)}\pi^{h^0(X,V)}\Big(\pi^{-s}\Big)^{\chi(X,V)}d\mu_{X,n}(d)-\int_{\mathcal M_{X,n}(d)}
\Big(\pi^{-s}\Big)^{\chi(X,V)}d\mu_{X,n}(d)\\
=&\sum_{i\geq 0}\pi^{i-sd+ns(g-1)}\mathrm{Vol}\Big({W_{X,n}^{=i}(d)}\Big)-\pi^{-sd+n(g-1)s}
\mathrm{Vol}\Big({\mathcal M_{X,n}(d)}\Big).\end{aligned}$$
Then
a simple version of  rank $n$ zeta function of $X$ is defined by
$$\widetilde\zeta_{X,n}(s):=\sum_{\alpha\geq 0}{\Large \int}_{\mathcal M_{X,n}(\alpha\cdot n)}^{\#}
\Big(\pi^{h^0(X,V)}-1\Big)\cdot \Big(\pi^{-s}\Big)^{\chi(X,V)}\cdot d\mu.$$
Tautologically, using the Riemann-Roch, the duality and the vanishing for semi-stable bundles,
we have the following
 
\begin{thm}
({\bf Zeta Facts})
(i) ({\bf Rationality}) $\widetilde\zeta_{X,n}(s)$ is rational  in $T=t^n=\pi^{-ns}$ with $t=\pi^{-s}$

\noindent
(ii) ({\bf Functional Equation})
$\widetilde\zeta_{X,n}(1-s)=\widetilde\zeta_{X,n}(s)$

\noindent
(iii) ({\bf Singularities}) $\widetilde\zeta_{X,n}(s)$ admits only two singularities, all simple poles, at $s=0,1$ 
$$\mathrm{Res}_{s=1}\widetilde\zeta_{X,n}(s)=\mathrm{Vol}\Big(\mathcal M_{X,n}(0)\Big).$$
\end{thm}

Note that by Witten's volume formula and our parabolic reduction conjecture, if this were the right zeta,
the special values of this rational function would be independent of $X$. Consequently, the functions themselves would all be the same.

Similar statement can be applied to $\widehat\zeta_{M,n}(s)$. For this we have to use a result of Ramanujan-Harder and Deninger ([D]), which claims that certain analytic functions are uniquely determined by their special values at almost all positive integers. This then suggests that our high rank zetas are independent of $M$, an assertion which should be understood in the framework  of the existence of projective flat connections in the theory of conformal blocks. 
   
\vskip 0.50cm
\noindent
{\bf REFERENCES}
\vskip 0.20cm 
\noindent
[ACGH]  E. Arbarello, M. Cornalba, P. Griffiths, \& J. Harris, {\it  Geometry of algebraic curves.} Vol. I. Grundlehren der Mathematischen Wissenschaften 267. Springer-Verlag, New York, 1985
\vskip 0.20cm 
\noindent
[Ad] K. Adachi, Parabolic Reduction, Stability and Volumes of Fundamental Domains: Examples of Low Ranks, Master Thesis, Kyushu Univ., 2012
\vskip 0.20cm 
\noindent
[AGBMNV] L. Alvarez-Gaume, J.-B. Bost, G. Moore, P. Nelson, \& C. Vafa,  
Bosonization on higher genus Riemann surfaces. Comm. Math. Phys. 112 (1987), no. 3, 503--552
\vskip 0.20cm 
\noindent
 [AB] M. Atiyah \& R. Bott, The Yang-Mills equations over Riemann surfaces. Philos. Trans. Roy. Soc. London Ser. A 308 (1983), no. 1505, 
 523-615
\vskip 0.20cm 
\noindent 
[A] J. Arthur, A trace formula for reductive groups. I. 
Terms associated to classes in $G({\mathbb Q})$. Duke Math. J. {\bf 45} (1978), 
no. 4, 911--952
\vskip 0.20cm 
\noindent
[BV] J.-M. Bismut \& E. Vasserot,  The asymptotics of the Ray-Singer analytic torsion associated with high powers of a positive line bundle. 
Comm. Math. Phys. 125 (1989), no. 2, 355--367
\vskip 0.20cm 
\noindent
[D] C. Deninger, How to recover an L-series from its values at almost all positive integers. Some remarks on a formula of Ramanujan. Proc. Indian Acad. Sci. Math. Sci. 110 (2000), no. 2, 121--132
\vskip 0.20cm
\noindent
[DR] U.V. Desale \& S. Ramanan, Poincare polynomials of the variety of stable bundles, Math. Ann 26 (1975) 233-244
\vskip 0.20cm
\noindent
[F]  G. Faltings, {\it  Lectures on the arithmetic Riemann-Roch theorem.}  Annals of Mathematics Studies, 127. Princeton University Press, Princeton, NJ, 1992
\vskip 0.20cm
\noindent
[F] J. Fay, {\it Kernel functions, analytic torsion, and moduli spaces} 
Mem. Amer. Math. Soc. 96 (1992), no. 464
\vskip 0.20cm
\noindent
[GS] G. van der Geer \& R.  Schoof, Effectivity of Arakelov divisors and the theta divisor of a number field. Selecta Math. (N.S.) 6 (2000), no. 4, 377--398
 \vskip 0.20cm
\noindent 
[HN] G. Harder \& M.S. Narasimhan, On the cohomology groups of moduli spaces of vector bundles on curves, Math. Ann. 212, 215-248 (1975)
 \vskip 0.20cm
\noindent
[KW] H.H. Kim \& L. Weng,  Volume of truncated fundamental domains. Proc. Amer. Math. Soc. 135 (2007), no. 6, 168--1688
\vskip 0.20cm
\noindent
[K] M. Kontsevich \& Y. Soibelman, Lecture Notes on Wall-crossing, Dec 2011
\vskip 0.20cm
\noindent
[L] L. Lafforgue, {\it Chtoucas de Drinfeld et conjecture de 
Ramanujan-Petersson}. Asterisque No. 243 (1997)
\vskip 0.20cm
\noindent
[LR] G. Laumon \& M. Rapoport, The Langlands lemma and the Betti numbers of stacks of G-bundles on a curve. Internat. J. Math. 7 (1996), no. 1, 29--45
\vskip 0.20cm
\noindent
[La] R. Langlands, The volume of the fundamental domain for some arithmetical 
subgroups of Chevalley groups, Proc. Sympos. Pure Math. 9, AMS (1966) pp.143--148
\vskip 0.20cm 
\noindent
[RR] S. Ramanan  \& A. Ramanathan, Some remarks on the instability flag. 
Tohoku Math. J. (2) 36 (1984), no. 2, 269--291
\vskip 0.20cm 
\noindent
[RS] D. Ray \& I. Singer,  Analytic torsion for complex manifolds. Ann. of Math. (2) 98 (1973), 154--177
vskip 0.20cm 
\noindent
[Weil] A. Weil, {\it Adeles and algebraic groups} Progress in Mathematics, 23. Birkh\"auser, Boston, Mass., 1982
\vskip 0.20cm 
\noindent
[W0] L. Weng, Geometry of Numbers,  arXiv:1102.1302 
\vskip 0.20cm 
\noindent
[W1] L. Weng, Non-abelian zeta functions for function fields. {\it Amer. J. Math}. 127 (2005), no. 5, 973--1017.
 \vskip 0.20cm 
\noindent 
[W2] L. Weng, A geometric approach to $L$-functions. {\it The Conference on $L$-Functions,} 219--370, World Sci. Publ., Hackensack, NJ, 2007.
 \vskip 0.20cm 
\noindent 
[W3] L. Weng, Symmetries and the Riemann Hypothesis, {\it Advanced Studies in Pure Mathematics 58}, Japan Math. Soc., 173-224, 2010
 \vskip 0.20cm 
\noindent 
[W4] L. Weng, Counting Bundles, arXiv:1202.0869
 \vskip 0.20cm 
\noindent 
[W5] L. Weng,  Zeta Functions for Elliptic Curves I: Counting Bundles, arXiv:1202.0870
 \vskip 0.20cm 
\noindent 
[W6] L. Weng, Zeta functions for function fields, arXir: 1202.3183
 \vskip 0.20cm 
\noindent 
[W7] L. Weng, Special uniformity of zeta functions I. Geometric aspect,  arXiv:1203.2305 
 \vskip 0.20cm 
\noindent 
[W8] L. Weng, Parabolic Reduction, Stability and The Mass, in preparation
 \vskip 0.20cm 
\noindent 
[WZ] L. Weng \& D. Zagier, High rank zeta functions for elliptic curves and the Riemann Hypothesis, in preparation
 \vskip 0.20cm 
\noindent 
[WZ2] L. Weng \& D. Zagier,  Parabolic reduction, stability and the mass: Vector bundles over curves, in preparation
 \vskip 0.20cm 
\noindent 
[Wi] E. Witten, On quantum gauge theories in two dimensions Comm. Math. Phys. Volume 141, Number 1 (1991), 153--209.
\vskip 0.20cm
\noindent
[Z] D. Zagier, Elementary aspects of the Verlinde formula and the Harder-Narasimhan-Atiyah-Bott formula, in {\it Proceedings of the Hirzebruch 65 Conference on Algebraic Geometry}, 445-462 (1996)

\vskip 4.50cm
\noindent
{\bf Lin WENG\footnote{We would like to thank Zagier for his kind arrangement, MPIM, Bonn for providing us with a nice and stimulating research environment,  Deninger and Hida for their constant encouragements, and in advance the reader who would tolerant our style of presentation. (After all, even the precision is always the heart of our culture, originality should be valued most.)

This work is partially supported by JSPS.
}},
 Institute for Fundamental Research, The $L$ Academy {\it and}

\noindent
Graduate School of Mathematics, Kyushu University, Fukuoka, 819-0395, 
JAPAN

\noindent
E-Mail: weng@math.kyushu-u.ac.jp

\end{document}